\documentclass[12pt]{article}

\usepackage{latexsym}
\usepackage{amsmath}
\usepackage{amsfonts}
\usepackage{amssymb}
\usepackage{graphics}

\newcommand{\bean}{\begin{eqnarray}}
\newcommand{\eean}{\end{eqnarray}}
\newcommand{\bea}{\begin{eqnarray*}}
\newcommand{\eea}{\end{eqnarray*}}
\newcommand{\bsa}{\begin{subarray}{c}}
\newcommand{\esa}{\end{subarray}}
\newcommand{\bi}{\begin{itemize}}
\newcommand{\ei}{\end{itemize}}

\newtheorem{lemma}{Lemma}
\newtheorem{thm}[lemma]{Theorem}

\newtheorem{cor}[lemma]{Corollary}

\begin{document}
\title{On the canonical decomposition of generalized modular functions}
\author{Geoffrey Mason\thanks{%
Supported by the NSF, NSA, and the Committee on Research at the University
of California, Santa Cruz} \\
Department of Mathematics, \\
University of California Santa Cruz, \\
CA 95064, U.S.A. \and Winfried Kohnen\\
Mathematisches Institut der Universit\"at Heidelberg\\
INF 288, D-69120 Heidelberg\\
Germany}
\date{}
\maketitle

\begin{abstract}
\noindent
 The authors have conjectured (\cite{KoM}) that if a normalized generalized modular function
(GMF) $f$, defined on a congruence subgroup $\Gamma$, has integral Fourier coefficients, then $f$ is classical in the sense that some power $f^m$ is a modular function on $\Gamma$. A strengthened form of this conjecture was proved (loc cit) in case the divisor of $f$ is \emph{empty}. In the present paper we study the canonical decomposition of a normalized parabolic GMF  $f = f_1f_0$ into a product of normalized parabolic GMFs $f_1, f_0$ such that $f_1$ has \emph{unitary character} and $f_0$ has \emph{empty divisor}. We show that the strengthened form of the conjecture holds if the first "few" Fourier coefficients of $f_1$ are algebraic. We deduce proofs of  several new cases of the conjecture, in particular if either $f_0=1$ or if the divisor of $f$ is concentrated at the cusps of $\Gamma$.
\end{abstract}

\section{Introduction}

Let $\Gamma\subset \Gamma_1:=SL_2({\bf Z})$ be a congruence subgroup and let $\mathcal{H}$ be the complex upper half-plane. We consider \emph{generalized modular functions} of weight zero (GMF's) on $\Gamma$. These are holomorphic functions $f:{\cal H}\rightarrow {\bf C}$ which satisfy
$$f(\gamma\circ z)=\chi(\gamma)f(z) \qquad (\forall \gamma\in \Gamma)$$
where $\chi:\Gamma\rightarrow \mathbb{C}^*$ is a (not necessarily unitary) character, and which are meromorphic at the cusps. We call $f$ a \emph{parabolic} generalized modular form (PGMF) if it also satisfies $\chi(\gamma)=1$ for all parabolic elements $\gamma\in \Gamma$ of trace 2.  If $f$ is a GMF then some power $f^m$ of $f$ is a PGMF. For further details we refer to \cite{KM1}. In this paper we deal mainly with PGMFs.

\medskip
At the cusp infinity, a PGMF $f$ has a Fourier expansion
$$f(z)=\sum_{n\geq h}a(n)q_N^n \qquad (0<|q_N|<\epsilon)$$
for appropriate $h\in \mathbb{Z}, N\in \mathbb{N}$  and where $q_N:= e^{2\pi iz/N}\,\, (z\in {\cal H})$.
\medskip
We shall call $f$ \emph{normalized} if $a(h)=1$.
\medskip
According to \cite{KM2}, each normalized PGMF $f$ on $\Gamma$ has a \emph{canonical decomposition}
$$f=f_1f_0 \leqno (1)$$
where $f_1$ and $f_0$ are normalized PGMF's on $\Gamma$, $f_1$ has \emph{unitary} character, and the divisor of $f_0$ is \emph{empty}. The canonical decomposition is indeed unique, which amounts to the assertion (loc cit) that a PGMF with empty divisor and unitary character 
is constant. 
Note that it follows from these conditions and our assumptions that
$$f_1=q_N^h+\dots, \qquad f_0=1+\dots .$$
\medskip
It was conjectured in \cite{KoM}, Sect. 1, that if a normalized PGMF $f$ on the Hecke congruence subgroup $\Gamma_0(N)$ of level $N$ has {\it integral} Fourier coefficients, then it must be classical, i.e. the character $\chi$ is of finite order. A proof of this conjecture would have some very important consequences in rational conformal field theory, as explained in \cite{KoM}. 

\medskip
In the present paper we shall show that the conjecture  is equivalent to requiring that the first ``few'' Fourier coefficients of the function $f_1$ in the decomposition (1) are algebraic numbers.  We will actually prove slightly stronger statements regarding the hypothesis on the Fourier coefficients of $f$, requiring only that they are rational and $p$-integral for almost all primes $p$. This result implies several new cases of the conjecture: if $\chi$ is unitary (i.e. $f_0 = 1$), or if the divisor $div(f)$ of $f$ is supported at the cusps of $\Gamma$. Some special cases of the second assertion, including the case when $div(f)$ is \emph{empty} (i.e. $f_1=1$), were established in \cite{KoM}. Our present results are valid for an arbitrary congruence subgroup $\Gamma$.

\medskip
Apart from the decomposition (1), there are two main
ingredients to the proof of our results: the first is a theorem of Scholl-Waldschmidt (\cite{Sc}, \cite{W}) on the transcendence of canonical differentials of the third kind on modular curves;  the second one is the result on PGMF's with empty divisors \cite{KoM}
already mentioned above,  whose proof largely depends on the analytic theory of Dirichlet series.

\medskip
The paper is organized as follows. We give the proof of the main results, Theorem \ref{thm1} and Corollary \ref{thm2}, in Section 3. In Section 4 we consider the action of complex conjugation on PGMFs and their characters. For example, we show that if the Fourier coefficients of a PGMF $f$ are real then so are those of $f_0$ and $f_1$.

\bigskip
\noindent {\it Acknowledgements.} The authors thank Henri Darmon and Jan H. Bruinier for useful conversations. The results of the present paper emerged from discussions  during the stimulating Workshop on noncongruence modular forms, organized by Ling Long and Winnie Li in August 2009,
at the American Institute of Mathematics in Palo Alto.

\section{Statement of results}
\bigskip 
We define 
$$\kappa_\Gamma:= \lbrack {1\over 6}[P\Gamma_1:P\Gamma]\rbrack +1-c_\Gamma.$$
Here $P\Gamma_1:=PSL_2({\bf Z})$ and $P\Gamma$ denotes the image of $\Gamma$ under the natural projection $\Gamma_1\rightarrow P\Gamma_1$. Furthermore
$[x]\, (x\in {\bf R})$ denotes the greatest integer function and $c_\Gamma$ is the number of cusps of $\Gamma$. \medskip
The main result of the paper is the following:
\begin{thm}\label{thm1} Let $f$ be a normalized PGMF on $\Gamma$ whose Fourier coefficients $a(n)$ are rational for all $n$ and are $p$-integral for all but a finite number of primes $p$. Let (1) be the canonical decomposition of $f$. Then the character $\chi$ of $f$ is of finite order if, and only if, the Fourier coefficients $a_1(n)$ of $f_1$ are algebraic for $h\leq n \leq\kappa_\Gamma +h$.
\end{thm}

\bigskip
Theorem \ref{thm1} has the following consequence.
\begin{cor}\label{thm2} Let $f$ be a PGMF on $\Gamma$ and let $\chi$ be the character of $f$.
Assume that the Fourier coefficients $a(n)$ are rational for all $n$ and are $p$-integral for all but a finite number of primes $p$. Then
$\chi$ has finite order if either of the following conditions hold:
\begin{eqnarray*}
&&(a) \ \mbox{$\chi$ is \emph{unitary}}, \\
&&(b) \ \mbox{The divisor of $f$ is concentrated at the cusps of $\Gamma$.}
\end{eqnarray*}
\end{cor}
We point out (\cite{KM1}) that the condition on the divisor in part (b) is equivalent to the assumption that the logarithmic derivative $f'/f$, which is generally a meromorphic modular form of weight $2$ on 
$\Gamma$, is in fact \emph{holomorphic}.

\section{Proof of Theorem \ref{thm1}}
\bigskip
Without loss of generality we may assume that $\Gamma=\Gamma(N)$ is the principle congruence subgroup of $\Gamma_1$ of level $N$. 

\medskip
In one direction the conclusion is easy. Indeed, assume that $\chi$ has finite order $m$. Then $f^m$ has trivial character, so that
$$f^m=f_1^m\cdot f_0^m=f_1^m\cdot 1$$
must be the canonical decomposition of $f^m$. Therefore $f$ and $f_1$ differ by an $m$-th root of unity, and
since they are both normalized then they are equal and all Fourier coefficients of $f_1$ are algebraic, indeed rational.

\medskip
Now supppose that $a_1(n)$ is algebraic for $h\leq n \leq\kappa +h$, where we have abbreviated $\kappa:=\kappa_\Gamma$. Solving recursively in (1) for the Fourier coefficients $a_0(n)$ of $f_0$, we see that our assumption implies that each $a_0(n) \ (0\leq n \leq\kappa)$ is also algebraic.

\medskip
Let
$${{2\pi i}\over N}g_0:={{f_0'}\over {f_0}}\leqno (2)$$
be the logarithmic derivative of $f_0$ and write
$$g_0=\sum_{n\geq 1}b_0(n)q_N^n.$$
Then we see that each $b_0(n)\, (1\leq n \leq\kappa)$ is algebraic. Note that $g_0$ is a cusp form of weight 2 on $\Gamma$ with trivial character, since $div(f_0)=\emptyset$ \cite{KM1}. 

\medskip
We now assert
\begin{lemma}\label{lemma1}All Fourier coefficients $b_0(n)\, (n\geq 1)$ are contained in a 
finite extension  $K/\mathbb{Q}$.
\end{lemma}
\begin{pf} Our claim is essentially well-known,  and follows from linear algebra combined with the valence formula and the fact that $\Gamma=\Gamma(N)$. However, for the reader's convenience we will give a detailed proof.

\medskip
First recall that the valence formula says that the sum of the orders (measured in the appropriate local variables) of a non-zero cusp form of weight 2 on $\Gamma$ on the complete modular curve $X_\Gamma:=\overline {\Gamma\backslash {\cal H}}$ is equal to $\lbrack {1\over 6}[P\Gamma_1:P\Gamma]\rbrack$.

\medskip
Since $\Gamma =\Gamma(N)$, it is well-known that the space $S_2(\Gamma)$ of cusp forms of weight 2 for $\Gamma$ has a basis $\lbrace g_1,\dots ,g_d\rbrace$ of functions with rational (in fact integral) Fourier coefficients (\cite{Sh}, Thm. 3.52). Note that the valence formula implies that $d\leq \kappa$.

\medskip
We write $g_0$ as a linear combination of the $g_\nu\, (1\leq\nu\leq d)$. Bearing in mind that the 
first $\kappa$ Fourier coefficients of $g_0$ are algebraic, our claim will follow if we can show that the 
$\kappa \times d$ matrix $A$ whose columns consists of the first $\kappa$ Fourier coefficients of $g_1, \dots , g_d$ has maximal rank.

\medskip
To show this we argue as follows. Let $\langle \, , \, \rangle$ denote the usual inner product 
on $S_2(\Gamma)$. Let $P_n\, (n\geq 1)$ be the $n$-th ``Poincar\'e series'' in $S_2(\Gamma)$ with
respect to $\langle \, , \,\rangle$, i.e. $P_n$ is the dual of the functional that sends a cusp form $g\in S_2(\Gamma)$ to its $n$-th Fourier coefficient $a_g(n)$. By the valence formula, $\lbrace P_1, \dots ,P_\kappa\rbrace$ generates $S_2(\Gamma)$. Hence,  there exists a basis $\lbrace P_{n_1},\dots ,P_{n_d}\rbrace$ where $1\leq n_\nu \leq \kappa$ for all $\nu$.

\medskip
On the other hand, let $\ell\in S_2^*(\Gamma)$ be any functional. Then by standard duality, there exists $L\in S_2(\Gamma)$ such that $\langle g,L\rangle =\ell(g)$ for all $g$. Writing $L$ in terms of the basis $\lbrace P_{n_1},\dots ,P_{n_d}\rbrace$, we see that $\ell$ is a linear combination of $a(n_1),\dots ,a(n_d)$, hence the latter functionals form a basis of $S_2^*(\Gamma)$.

\medskip
From the above it follows that the $\kappa \times d$ matrix $B$ whose columns are the first $\kappa$ Fourier coefficients of $P_{n_1},\dots ,P_{n_d}$ has maximal rank (the rows with indices $n_1, \dots ,n_d$ are linearly independent). Hence the same is true for $A$, since $A$ is obtained from $B$ by multiplying with an invertible $d \times d$ matrix.
This completes the proof of the Lemma. $\hfill \Box$
\end{pf}

\bigskip
In (2) we now solve recursively for the $a_0(n)$. Using Lemma 1, we see that $a_0(n)\in K$ for all $n$. Therefore by (1),  each $a_1(n)$ also lies in $K$.

\medskip
We put
$${{2\pi i}\over N}g_1:={{f_1'}\over {f_1}}. \leqno (3)$$
Then $g_1$ also has Fourier coefficents in $K$, and $g_1$ is a meromorphic modular form of weight 2 on $\Gamma$ with trivial character.  It has at worst simple poles in ${\cal H}$ with integral residues, and is holomorphic at the cusps \cite{KM1}.

\medskip
We let
$$D:=div(f).$$
Then $deg D=0$ \cite{KM1}, while $div(f)=div(f_1)$ by hypothesis.

\medskip
The form $g_1$ defined by (3) gives rise to an abelian differential of the third kind
$$\omega_1:={{2\pi i}\over N}g_1 dz$$
on $X_\Gamma$ with residue divisor $D$.

\begin{lemma}\label{lemma2} The divisor $D$ is defined over a number field.
\end{lemma}
\begin{pf} This essentially is well-known: indeed, the Galois group operates on meromorphic differentials and this operation is compatible with the formation of residue divisors. However, for the convenience of the reader we again give a detailed proof.
We shall prove somewhat more, namely that each point $P$ in the support of $D$ is already fixed by 
$Gal(\overline{ \mathbb{Q}}/L)$ where $L/\mathbb{Q}$ is a finite extension.

\medskip
The modular curve $X_\Gamma$ is defined over a number field $F$ (in fact, over $\mathbb{Q}(\zeta_N)$ where $\zeta_N=e^{2\pi i/N}$). Since the cusps of $\Gamma$ are defined over a number field, we may suppose that $P$ is contained in the ``open part'' $Y_{\Gamma}:= X_{\Gamma}\setminus \lbrace cusps\rbrace$ of the modular curve.

 \medskip
The function field $F_{X_\Gamma}$ of $X_\Gamma/F$ is a finite extension of $F(j)$ where $j$ is the classical modular invariant. By the ``theorem of the primitive element''  there exists a modular function $t$ for $\Gamma$ such that $F_{X_\Gamma} =F(j,t)$, and $t$ satisfies an algebraic equation over $F(j)$. The points of $Y_{\Gamma}$ then can be parametrized as $(j(z),t(z))\, (z\in \Gamma\backslash {\cal H})$.

\medskip
Suppose that $P$ ``corresponds'' to $(j(z),t(z))\, (z\in {\cal H})$. It is sufficient to show that $j(z)$ is algebraic over $\mathbb{Q}$. Because $P$ is not a cusp then $z$ is necessarily a zero of $f$, and hence a simple pole of $g_1$.
Let
$$G_1:= \frac{g_1^6}{\Delta},$$
where $\Delta$ is the usual discriminant function of weight 12 on $\Gamma_1$. Then $G_1$ is a modular function on $\Gamma$ with a pole at $z$. 
Because $g_1$ has Fourier coefficients in $K$, the same is true of $G_1$.
Observe that each of the translates
$$(G_1|\gamma)(z):=G_1(\frac{az+b}{cz+d}), \qquad \gamma= \left(\begin{array}{cc}a & b \\c & d\end{array}\right)\in \Gamma_1,$$
also have Fourier coefficients in a number field. Once again, this is a standard result. The proof is
based on  the ``$q$-expansion principle" for congruence subgroups (\cite{DR}, VII, 4.8). (A more general proof that works for noncongruence subgroups can be found in \cite{KL}, Proposition A1.)

\medskip
We may, and shall, choose $c \in \mathbb{Q}$ such that \emph{none} of the modular functions
$(G_1|\gamma)(z) - c$ have a zero at $z$. Now consider the ``norm''
$$\prod_{\gamma\in P\Gamma\backslash P \Gamma_1} (G_1 - c)|\gamma.$$
It is a modular function on $\Gamma_1$ with coefficients in $K$, and hence is a rational function 
$A(j)/B(j)$ with $A(j), B(j) \in K[j]$. By construction there is a pole at $z$, so that $B(j(z))=0$. The algebraicity of $j(z)$ follows, and the proof of the Lemma is complete.  $\hfill \Box$
\end{pf}

\bigskip
Next, observe that $\omega_1$ is a {\it canonical} differential on $X_\Gamma$, i.e.
$$\Re(\int_\sigma\omega_1)=0$$
for all $\sigma\in H_1(U;{\bf Z})$ where $U:=X_\Gamma\setminus supp D$. Indeed, this is equivalent to our assumption that $f_1$ has unitary character, bearing in mind the relation (3).
We may now invoke an important theorem of Scholl-Waldschmidt \cite{Sc}, \cite{W}: since $g_1$ has Fourier coefficients in a number field, $D$ has finite order in the divisor class group. Thus there exists $m\in {\mathbb{N}}$ and a modular function $h_1$ on $\Gamma$ such that
$$div({{f_1^m}\over {h_1}})=\emptyset.$$
Therefore $f_1^m=h_1$, since $f_1$ and $h_1$ both have unitary character.

\medskip
We may normalize $h_1$ to have Fourier coefficients in the compositum of ${\bf Q}(\zeta_N)$ and the field of definition of the algebraic divisor $D$ (cf. Lemma \ref{lemma2}) and to have leading non-zero term equal to 1.
From (1) we now obtain
$$f^m=h_1\cdot h_0, \leqno (4)$$
where 
$$h_0:=f_0^m.$$
Since $div(h_0)=\emptyset$ and $h_1$ has trivial character, it follows that (4) is the canonical decomposition of $f^m$.
\begin{lemma}\label{lemmahconj} Let $h=\sum b(n)q_N^n$ be a modular function on $\Gamma =\Gamma(N)$ with Fourier coefficients in a number field, and let $\sigma \in Gal(\overline {\mathbb{Q}}/\mathbb{Q})$. Then $h^\sigma:=\sum b(n)^\sigma q_N^n$ is a modular function on $\Gamma$.
\end{lemma}
\begin{pf} This is well-known. Indeed, arguing as in the proof of Lemma \ref{lemma2},  we see that $j(z)$ is algebraic if $z\in {\cal H}$ is a pole of $h$. We therefore conclude that there exists a polynomial $P$ with algebraic coefficients and a large positive integer $M$ such that
$$H:=\Delta^MP(j)h$$
is a cusp form of weight $k\geq 2$ on $\Gamma$. The latter space has a basis of functions with rational Fourier coefficients \cite{Sh}, so that $H$ is a $\overline{\mathbb {Q}}$-rational linear combination of such functions (an argument similar to that used in proof of Lemma \ref{lemma1} with $2$ replaced by $k$ is valid). The result follows from this. $\hfill \Box$
\end{pf}

\bigskip
All coefficients of the individual functions on both sides of (4) are contained in a finite Galois extension $F/\mathbb{Q}$. In (4) we now take ``norms'' (products of Galois conjugates of 
$\sigma\in Gal(F/\mathbb{Q}))$ to obtain an equation
$$f^{mr}=H_1\cdot H_0, \leqno (5)$$
where $r=[F: \mathbb{Q}]$, $H_1$ and $H_0$ are normalized and have rational Fourier coefficients,  and (thanks to Lemma \ref{lemmahconj})  $H_1$ has trivial character.

\medskip
Now notice that $div(H_0)=\emptyset$. Indeed, if ${{2\pi i}\over N}w \in S_2(\Gamma)$ corresponds to the normalized PGMF $v$ upon taking logarithmic derivatives, then the cusp form ${{2\pi i}\over N}w^\sigma$ corresponds to $v^\sigma$ (where of course the action of Galois elements $\sigma$ is defined in the same way as above).

\medskip
Since $H_1$ has rational coefficients and $\Gamma =\Gamma(N)$, the coefficients of $H_1$ must in fact be $p$-integral for almost all primes $p$ (\cite{Sh}). Since $H_1$ is normalized, the same holds for $H_1^{-1}$ and therefore also for $H_0$ by (5).
So we have arrived at the situation that $H_0$ is a normalized PGMF with \emph{rational Fourier coefficients} and \emph{empty divisor}. By \cite{KoM}, Theorem 2,  we conclude that $H_0=1$. Thus $f^{mr}=H_1$ has trivial character, i.e. the character of $f$ is of finite order. This concludes the proof of Theorem \ref{thm1}.

\bigskip
We turn to the proof of  Corollary \ref{thm2}. Suppose first that $\chi$ is \emph{unitary}. Then $f_0=1$, so that
$f=f_1$ has rational Fourier coefficients by hypothesis. Therefore, $\chi$ has finite order by Theorem \ref{thm1}. This proves part (a) of the Corollary. As for part (b), because the divisor 
$D = div(f)=div(f_1)$ is assumed to be concentrated at the cusps, $D$ is defined over a number field and the Manin-Drinfeld theorem tells us that $div(f)$ has finite order in the divisor class group. Then
$f_1^m$ has algebraic Fourier coefficients for some integer $m$, hence $f_1$ does too. Now Theorem \ref{thm1} again tells us that $\chi$ has finite order. This completes the proof of part (b) of the Corollary.

\section{PGMFs with real Fourier coefficients}
In this Section we briefly consider the action of complex conjugation on PGMFs.
Recall Hecke's  operator $K$, defined on  holomorphic functions in $\mathcal{H}$
 as follows (\cite{R}, Section 8.6)
 \begin{eqnarray*}
f|K (z) = \overline{ f(- \overline{z})}.
\end{eqnarray*}
 If $f$ is a PGMF on $\Gamma$ then so is $f|K$ (loc. cit.), and the $q$-expansions at the infinite cusp are related as follows:
$$
f(z) = \sum a(n)q_N^n, \ \ \ \ f|K(z)= \sum \overline{a(n)} q_N^n. \leqno(6)
$$

\bigskip
Set
\begin{eqnarray*}
J = \left(\begin{array}{cc}-1 & 0 \\0 & 1\end{array}\right).
\end{eqnarray*}
Note that
\begin{eqnarray*}
J\left(\begin{array}{cc}a & b \\ c & d\end{array}\right)J^{-1} &=& \left(\begin{array}{rr}a & -b \\ -c& d\end{array}\right).
\end{eqnarray*}

\medskip
In the following we will need to assume that $J$ normalizes $\Gamma$. From the last display, we see that  this holds, for many congruence subgroups, e.g., $\Gamma = \Gamma(N)$, or $\Gamma_0(N)$. We write $\gamma^J = J\gamma J^{-1}$ for 
$\gamma \in \Gamma$. A character $\chi$ of
$\Gamma$ may be `twisted' by $J$ to yield a second character $\chi^J$
defined by
\begin{eqnarray*}
\chi^J(\gamma) = \chi( \gamma^J).
\end{eqnarray*}
\begin{lemma}\label{lemmaR} Assume that $J$ normalizes $\Gamma$, and suppose that the PGMF $f$ is associated with the character $\chi$. 
Then $f|K$ is associated with the character $\bar{\chi}^J$. In particular, $f$ has real Fourier coefficients
if, and only if, $\chi = \bar{\chi}^J$.
\end{lemma}
\begin{pf} For $\gamma \in \Gamma$ we have
\begin{eqnarray*}
f|K(\gamma z) &=&  \overline{ f(- \overline{\gamma z})} =  \overline{ f(\gamma^J( -\overline{z}))}  = \bar{\chi}^J(\gamma)f|K(z).
\end{eqnarray*}
This proves the first assertion. The second follows from (6). $\hfill \Box$
\end{pf}

\bigskip
Suppose that $f$ is a PGMF with canonical decomposition (1),  and that $\chi_j$ is the character associated to $f_j, \ j=0, 1$. We have
$$ 
f|K = (f_1|K)( f_0|K). \leqno{(7)}
$$
By Lemma \ref{lemmaR}, $f_j|K$ has associated character $\overline{\chi_j}^J$. Because $\chi_1$ is unitary then so too is $\overline{\chi_1}^J$. Moreover, it is easy to see that $f_0|K$ has empty divisor. It follows from these comments that (7) is the canonical decomposition of $f|K$.

\medskip
If now we assume that $f$ has real Fourier coefficients, then $\chi = \overline{\chi}^J$ by Lemma
\ref{lemmaR}. By the uniqueness of the canonical decomposition, we can conclude that 
$\overline{\chi_j}^J = \chi_j$ for $j =0, 1$. Applying Lemma \ref{lemmaR} once more, we arrive at
\begin{lemma}\label{realdecomp}
Suppose that $f$ has real Fourier coefficients. Then $f_1$ and $f_0$ have real Fourier coefficients.
$\hfill \Box$
\end{lemma} 

\medskip
The condition $\overline{\chi_j}^J = \chi_j$ places strong restraints on the characters $\chi_j$.
For example, in the case of the unitary character $\chi_1$ we have 
$\chi_1(\gamma^J \gamma)=1$ for $\gamma \in \Gamma.$ 

\medskip
Finally, notice that Lemma \ref{realdecomp} applies in the context of Theorem \ref{thm1}. If we could
replace the real field  by a number field in the statement of the Lemma, the main conjecture would follow.

\bigskip
\noindent
e-mail addresses of authors:\\
\noindent {\it Winfried Kohnen}  winfried@mathi.uni-heidelberg.de\\
\noindent {\it Geoffrey Mason} gem@cats.ucsc.edu

\end{document}